\documentclass[reqno]{amsart}
\usepackage{amssymb,latexsym,epsf}

\newtheorem{thm}{Theorem}[section]

\makeatletter
\@addtoreset{figure}{section}
\def\thefigure{\thesection.\@arabic\c@figure}
\def\fps@figure{h, t}
\@addtoreset{table}{bsection}
\def\thetable{\thesection.\@arabic\c@table}
\def\fps@table{h, t}
\@addtoreset{equation}{section}

\makeatother

\allowdisplaybreaks
\def\intprod{\mathbin{\hbox to 6pt{%
                 \vrule height0.4pt width5pt depth0pt
                 \kern-.4pt
                 \vrule height6pt width0.4pt depth0pt\hss}}}

\begin{document}

\title{On the stability of periodic 2$D$ Euler-$\alpha$ flows}

\author[S. Pekarsky]{Sergey Pekarsky}
\address{ {CDS\\California Institute of Technology, 107-81\\
Pasadena, CA 91125}}
\email{{sergey@cds.caltech.edu}}

\author[S. Shkoller]{Steve Shkoller}
\address{Dept. of Mathematics, University of California at Davis
\\  Davis, CA 95616-8633 }
\email{{shkoller@cds.caltech.edu}}


\date{September 1998; current version February 25, 1999}

\begin{abstract}
Sectional curvature of the group $\mathcal{D}_\mu (M)$
of volume-preserving 
diffeomorphisms of a two-torus with the $H^1$ metric is analyzed.
An explicit expression is obtained for the sectional
curvature in the plane spanned by two stationary flows,
$\cos (k, x)$ and $\cos (l, x)$.
It is shown that for certain values of the wave vectors $k$ and $l$
the curvature becomes positive for $\alpha > \alpha_0$,
where $0 < \alpha_0 < 1$ is of the order $1/k$.
This suggests that the flow corresponding to such geodesics becomes 
more stable as one goes from usual Eulerian description to the 
Euler-$\alpha$ model. 
\end{abstract}

\maketitle

\tableofcontents

\section{Introduction}

In Lagrangian mechanics a motion of a natural mechanical system
is a geodesic line on a manifold - configuration space in the
metric given by the difference of kinetic and potential energy.
The configuration space for the fluid motion in a domain $M$ is 
the  group $\mathcal{D}_\mu (M)$ of volume-preserving 
diffeomorphisms of $M$. The corresponding (Lie) algebra is the
algebra of divergence-free vector fields on $M$ vanishing on
the boundary. The standard (Euler) model of an ideal fluid corresponds 
to the kinetic energy being given by the $L^2$ norm of the fluid 
velocity on $M$. That is, the right-invariant metric on 
$\mathcal{D}_\mu (M)$ is defined in the following way:
its value at the identity of the group on a divergence-free vector
field $v$ from the algebra is given by 
$\langle v, v \rangle = \| v \|_{l^2} = \int_M (v,v) d x$.

Recently, a number of papers (see, e.g., \cite{HMR, S98, S99}) introduced
the so called averaged Euler equations for ideal 
incompressible flow on a manifold $M$.
The averaged Euler equations involve a parameter $\alpha$;
one interpretation is that they are obtained by temporally averaging
the Euler equations in Lagrangian representation over rapid
fluctuations whose amplitudes are of order $\alpha$.
The particle flows associated with these equations can be shown
to be geodesics on a suitable group of volume-preserving
diffeomorphisms but with respect to a right invariant $H^1$
metric instead of the $L^2$ metric.

The case of area-preserving diffeomorphisms of the two-dimensional
torus with a right invariant $L^2$ metric was analyzed by Arnold who 
showed (see, e.g. \cite{A66, AK}) that 
``in many directions the sectional curvature is negative''.
In this paper we consider geodesic stability problem for 
the group $\mathcal{D}_\mu (T^2)$ with a right invariant $H^1$
metric which is related to the average Euler flows.

The instability discussed in this paper is the exponential 
\emph{Lagrangian} instability of the motion of the fluid, not of its
velocity field. A stationary flow can be a Lyapunov stable
solution of Euler equations, while the corresponding motion of the
fluid is exponentially unstable. The reason is that a small
perturbation of the fluid velocity field can induce exponential
divergence of fluid particles.

\section{Instability of the Euler flow on $T^2$}

Here we review Arnold's results for the group $\mathcal{D}_\mu (T^2)$ 
with a right invariant $L^2$ metric closely following \cite{AK}.
Recall some standard notations.
Let $B$ denote the bilinear form on a Lie algebra ${\mathfrak g}$ 
defined by the relation 
$\langle B(\xi, \eta), \zeta \rangle = \langle \xi, 
[\eta, \zeta ] \rangle$, where $\xi, \eta, \zeta \in {\mathfrak g}$
$[\cdot, \cdot ]$ is the commutator in ${\mathfrak g}$ and 
$\langle \cdot, \cdot \rangle$ is the inner product
in the space ${\mathfrak g}$.

The (Riemannian) \emph{curvature tensor}
$R$ describes the infinitesimal transformation on a tangent space
obtained by parallel translation around an infinitely small
parallelogram. For $u, v, w \in T_{x_0} M$, the action of 
$R (u,v)$ on $w$ can be expressed in terms of covariant
differentiation as follows
\begin{equation}
  \label{cur_def}
R(u,v) w = (-\nabla_{\bar{u}}\nabla_{\bar{v}} \bar{w}
+\nabla_{\bar{u}}\nabla_{\bar{v}} \bar{w} +
+\nabla_{\{\bar{u},\bar{u}\}} \bar{w} ) |_{x=x_0},
\end{equation}
where $\bar{u},\bar{v},\bar{w}$ are any fields whose values at 
the point $x_0$ are $u,v,w$. 

The \emph{sectional curvature} of $M$ in the direction of the
two-plane spanned by any two vectors $u,v \in T_{x_0} M$ is
the value
\begin{equation}
  \label{sc_def}
  C_{u v} = \frac{\langle R(u,v)u, v \rangle}
{\langle u, u \rangle \langle v, v \rangle - \langle u, v \rangle^2}.
\end{equation}

Theorem $3.2$ of \cite{AK} gives explicit formulas for the
inner product, commutator, operation $B$, connection, and curvature of the
right invariant $L^2$ metric on the group $\mathcal{D}_\mu (T^2)$.
These formulas allow one to calculate the sectional curvature
in any two-dimensional direction.

The divergence-free vector fields that constitute the Lie algebra
of the group $\mathcal{D}_\mu (T^2)$ can be described by their
stream (Hamiltonian) functions with zero mean 
(i.e., $v = -\dfrac{\partial H}{\partial y} \dfrac{\partial}{\partial x} 
+\dfrac{\partial H}{\partial x} \dfrac{\partial}{\partial y}$).
Thus, the Lie algebra can be identified with the space of real
functions on the torus having zero average value \cite{AK}.
It is convenient to define such functions by their Fourier
coefficients and to carry out all calculations over ${\mathbb C}$.

Complexifying the Lie algebra one constructs a basis of this vector
space using the functions $e_k$ (where $k$, called a 
\emph{wave vector}, is a point of ${\mathbb R}^2$) whose value at a 
point $x$ of our complex plane is equal to $e^{\imath (k,x)}$.
This determines a function on the torus if the inner product 
$(k,x)$ is a multiple of $2 \pi$ for all $x \in \Gamma$. 
All such vectors $k$ belong to a lattice $\Gamma^*$ in ${\mathbb R}^2$,
and the functions $\{e_k | k \in \Gamma^*, k \ne 0 \}$ form a
basis of the complexified Lie algebra.

Consider the parallel sinusoidal steady flow
given by the stream function $\xi = \cos (k, x)$ 
and let $\eta$ be any other vector of the algebra, i.e. 
$\eta = \sum x_l e_l$, where $x_{-l} = \bar{x}_l$.
Theorem $3.4$ of \cite{AK} states that the curvature of the group 
${\mathcal D}_\mu (T^2)$ in any two-dimensional plane containing
the direction $\xi$ is \emph{non-positive} and is given by
\begin{equation}
  \label{Arn3.2}
  C_{\xi \eta} = \frac{S}{4} \sum_l a_{k l}^2 |x_l + x_{l+2k} |^2,
\end{equation}
where  
$a_{k l} = \dfrac{(k \times l)^2}{| k + l |}$, 
$k \times l = k_1 l_2 - k_2 l_1$ is the (oriented) area of the 
parallelogram spanned by
$k$ and $l$, and $S$ is the area of the torus.
Then, a corollary of this theorem states that the curvature
in the plane defined by the stream functions $\xi = \cos (k,x)$
and $\eta = \cos (l,x)$ is 
\begin{equation}
  \label{Arn3.3}
  C_{\xi \eta} = - (k^2 + l^2) \sin^2 \beta \sin^2 \gamma / 4S,
\end{equation}
where $\beta$ is the angle between $k$ and $l$, and $\gamma$ is
the angle between $k+l$ and $k-l$.

\section{Stable directions for the Euler-$\alpha$ flow on $T^2$}

In this section we present new results on the sectional curvature
of the group of area-preserving diffeomorphisms of a two-torus with
a right invariant $H^1$ metric in view of the application to
the Lagrangian stability analysis following Arnold \cite{A66}.
The foundations for these results were established in \cite{S98}
where the continuous differentiability of the geodesic spray
of $H^1$ metric on ${\mathcal D}^s_\mu(M)$ for an arbitrary Riemannian
manifold $M$ was proved.

We start with an analog of Theorem $3.2$ of \cite{AK}.
Define an operator 
$A^\alpha : {\mathbb R}^2  \rightarrow {\mathbb R}_+, \quad 
k \mapsto k^2 (1+\alpha^2 k^2)$. It corresponds to the $H^1$
norm in the Fourier space and is simply given by $k^2$ in the
case $\alpha = 0$ when the $H^1$ metric effectively becomes 
the $L^2$ metric.

\begin{thm}
The explicit formulas for the inner product, commutator, 
operation $B$, and connection of the
right invariant $H^1$ metric on the group $\mathcal{D}_\mu (T^2)$
have the following form:
\begin{equation}
  \label{thm_1}
  \langle e_k, e_l \rangle = A^\alpha (k) \delta_{k, -l}  
\end{equation}
\begin{equation}
  \label{thm_2}
[e_k, e_l] = (k \times l) e_{k+l} 
\end{equation}
\begin{equation}
  \label{thm_3}
B(e_k, e_l) = b_{k,l} e_{k+l}, \quad \operatorname{where}
\quad b_{k,l} = (k \times l) \frac{A^\alpha (k)}{A^\alpha (k+l)} 
\end{equation}
\begin{equation}
  \label{thm_4}
\nabla_{e_k} e_l = d_{k,k+l} e_{k+l}, \quad \operatorname{where}
\quad d_{k,k+l} = \frac{k \times l}{s} 
\left(1 - \frac{A^\alpha (k)-A^\alpha (l)}{A^\alpha (k+l)} \right).
\end{equation}
\end{thm}

Using the definition of the curvature tensor (\ref{cur_def})
we obtain
\begin{multline}
  \label{cur_coef}
  R_{k,l,m,n} \equiv \langle R(e_k, e_l) e_m, e_n \rangle  
= (-d_{l+m, k+l+m} d_{m, l+m} \\ + d_{k+m, k+l+m} d_{m, k+m}
+(k\times l) d_{m,k+l+m} ) A^\alpha (k+l+m) S.
\end{multline}

We do not write here the explicit expression for $R_{k,l,m,n}$ as
it is rather involved, but we note that it is non-zero only
in the case $k+l+m+n = 0$.
We analyze  a special case of the curvature
in the plane defined by the stream functions $\xi = \cos (k,x)$
and $\eta = \cos (l,x)$ (notice that the corresponding flow
is a solution of the averaged Euler equations).
Then the sectional curvature is determined only by two terms
(we ignore the scaling factor of the denominator in the definition 
(\ref{sc_def})):
$$
C^{H^1}_{\xi \eta} = \frac{1}{8} (R_{k,l,-k,-l}+R_{-k,l,k,-l})
$$
The computation gives an explicit formula
\begin{multline}
  \nonumber
  C^{H^1}_{\xi \eta} = \frac{S}{36} (k \times l)^2 
\left( 4 A^\alpha (k)+4 A^\alpha (l)- 3 A^\alpha (k+l)-3 A^\alpha (k-l)
\right. \\ \left.+\frac{(A^\alpha (k)-A^\alpha (l))^2}{A^\alpha (k-l)} 
+\frac{(A^\alpha (k)-A^\alpha (l))^2}{A^\alpha (k+l)}\right)
\end{multline}
which we rewrite in the following form
\begin{multline}
  \label{sc_res}
  C^{H^1}_{\xi \eta} = \rho^2 \{ A^\alpha (k+l) A^\alpha (k-l)
(4 A^\alpha (k)+4 A^\alpha (l)- 3 A^\alpha (k+l)-3 A^\alpha (k-l)) \\
+(A^\alpha (k)-A^\alpha (l))^2 (A^\alpha (k+l) + A^\alpha (k-l)) \},
\end{multline}
where $\rho^2 = \dfrac{S (k \times l)^2}
{36 A^\alpha (k+l) A^\alpha (k-l)}$ is a function of $k,l,\alpha$ and
is strictly positive. Hence, the sign of the curvature is determined
by the expression in the bracket, which is a cubic polynomial in
$\alpha^2$:
\begin{equation}
  \label{poly}
  B(\alpha, k,l) \equiv b_0 + b_1 \alpha^2 + b_2 (\alpha^2)^2 + 
  b_3 (\alpha^2)^3, 
\end{equation}
so that $C^{H^1}_{\xi \eta} = \rho^2 B(\alpha, k,l)$.

\begin{figure}
  \begin{center}
    \leavevmode
  \epsfxsize=3in
  \epsffile{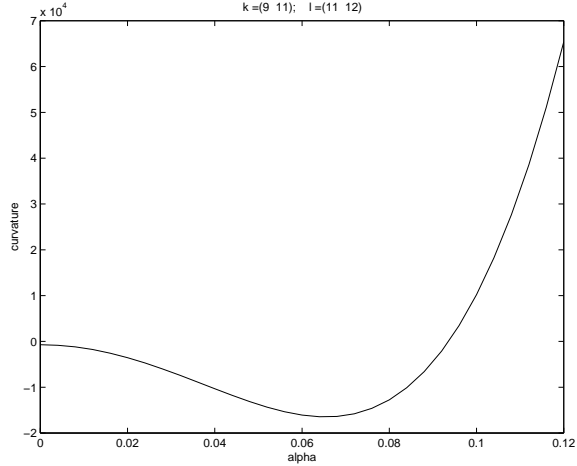}
    \caption{Sectional curvature (\ref{sc_res}) as a function of
     $\alpha$ for the case $k=(9,11), l=(11,12)$.}
    \label{fig:curv}
  \end{center}
\end{figure}

Numerical analysis of this complicated expression shows that
the sectional curvature becomes positive for some values of
$\alpha > \alpha_0$ when $k-l$ is small.
Fig. (\ref{fig:curv}) is representative of a typical behavior 
of the curvature as a function of $\alpha$ for $l = k + \epsilon$, 
where $\epsilon \ll k$ is small.
Based on this numerical evidence we analyze further analytically 
the case $l = k + \epsilon$, where $\epsilon \ll k$ is small.
Compute the coefficients $b_n$ in (\ref{poly}) as power
series in $\epsilon$
\begin{eqnarray}
  \label{Bcoef}
&  b_0  = -64 k^4 \epsilon^2 + 16 k^2 (k,\epsilon)^2 + 
  {\mathcal O}(\epsilon^4) \\
&  b_1  = -224 k^6 \epsilon^2 + 128 k^4 (k,\epsilon)^2 +
  {\mathcal O}(\epsilon^4) \\
&  b_2 = -640 k^8 \epsilon^2 + 320 k^6 (k,\epsilon)^2 +
  {\mathcal O}(\epsilon^4) \\
&  b_3 = 256 k^8 (k,\epsilon)^2 +
  {\mathcal O}(\epsilon^4) 
\end{eqnarray}

Notice that the coefficient of the highest degree is positive
while all the rest are negative. Hence, for $k > 1/\alpha$ it 
defines the leading term which increases with $\alpha$,
while the other coefficients are responsible for initial decrease
seen in Fig. (\ref{fig:curv}).
We summarize our result in the following theorem.

\begin{thm}
Consider the sectional curvature of the group $\mathcal{D}_\mu (T^2)$
equipped with the right invariant $H^1$ metric
in the plane defined by the stream functions $\xi = \cos (k,x)$
and $\eta = \cos (l,x)$, where $l=k+\epsilon$.
Then, for $|\epsilon|$ sufficiently small, for any $k$ there is an
$0 < \alpha_0(k) < 1$, such that for all $\alpha > \alpha_0(k)$ 
the corresponding sectional curvature is positive.
\end{thm}


\section*{Acknowledgments}

The authors would like to thank Jerrold E. Marsden for helpful
comments and the Center for Nonlinear Science of Los Alamos
National laboratory for providing a valuable setting where much 
of this work was performed.



\begin{thebibliography}{1234567}

\bibitem[A 66]{A66}{\sc V.I. Arnold}, 
Sur la geometrie differentielle des groupes de Lie de 
dimension infinie et ses applications a l'hydrodynamique des 
fluides parfaits, {\it Ann. Inst. Fourier}, Grenoble {\bf 16},
(1966), 319=361.

\bibitem[AK 98]{AK}{\sc V.I. Arnold and B. Khesin}, {\it
Topological Methods in Hydrodynamics}, Springer Verlag, New York,
1998.

\bibitem[HMR]{HMR}{\sc D.D. Holm, J.E. Marsden, and T.S. Ratiu},
Euler-Poincar\'{e} equations and semidirect products with application
to continuum theories, {\em Adv. in Math.} {\bf 137}, (1998), 1-81.

\bibitem[L 88]{L88}{\sc A.M. Lukatski}, Structure of the
curvature tensor of the group of measure-preserving diffeomorphisms of
a compact to-dimensional manifold, {\it Sibirskii Math. J.} {\bf 29}(6)
(1988), 95-99. 

\bibitem[S 98]{S98}{\sc S. Shkoller}, Geometry and curvature 
of diffeomorphism groups with $H^1$ metric and mean hydrodynamics,
{\it J. Functional Analysis}, {\bf 160}, (1998) 337-365

\bibitem[S 99]{S99}{\sc S. Shkoller}, The geometry and analysis
of non-Newtonian fluids and vortex methods, {\it preprint} (1999)


\end{thebibliography}
\end{document}